\documentclass[%
 preprint,
 amsmath,amssymb,
 aps,
]{revtex4-1}

\usepackage{algorithm}
\usepackage[]{algorithmicx}
\usepackage{algpseudocode}
\usepackage{graphicx}

\begin{document}


\title[Bounding the Radius of Convergence of Analytic Functions]{Bounding the Radius of Convergence of Analytic Functions} 
\author{Adam S. Jermyn}
\affiliation{Institute of Astronomy, University of Cambridge, Cambridge CB3 0HA, UK.}
\date{\today}

\begin{abstract}

Contour integration is a crucial technique in many numeric methods of interest in physics ranging from differentiation to evaluating functions of matrices.
It is often important to determine whether a given contour contains any poles or branch cuts, either to make use of these features or to avoid them.
A special case of this problem is that of determining or bounding the radius of convergence of a function, as this provides a known circle around a point in which a function remains analytic.
We describe a method for determining whether or not a circular contour of a complex-analytic function contains any poles.
We then build on this to produce a robust method for bounding the radius of convergence of a complex-analytic function.

\end{abstract}

\maketitle

\section{Introduction}

Contour integration is used extensively in numerical analysis.
In many cases however it is necessary to know something about the pole structure of a function.
This is why, for instance, many contour eigenvalue methods require that all poles lie on the real axis\ \citep{sakurai2007}.
Newer methods have been developed in some domains without these requirements\ \citep{7464870}, but these are the exception rather than the rule.

A special case of this problem is that of choosing a contour in which the function is analytic.
This arises in numerical differentiation via contour integration, which requires that the contour contain no poles\ \citep{Fornberg:1981:NDA:355972.355979}.
Similarly being able to exclude regions from eigenvalue searches can aide with the convergence of these methods\ \citep{7464870}.

The radius of convergence of a function $f(z)$ about the point $z_0$ is the radius $r(z_0) = |z - z_0|$ of the largest circle in the complex plane such that
\begin{equation}
	f(z) = \sum_{n=0}^\infty a_n (z - z_0)^n
	\label{eq:series}
\end{equation}
converges.
Equivalently,
\begin{equation}
	r(z_0) = |z_1 - z_0|,
\end{equation}
where $z_1$ is the nearest singularity of $f(z)$ to $z_0$.
The latter definition is more helpful, as computing power series numerically to large order is expensive, while searching for singularities is potentially more straightforward.

In this work we present a simple and robust method for determining whether or not a given circular contour contains any poles.
In Section\ \ref{sec:integrate} we show that the question of whether or not a circular contour contains any poles is equivalent to a statement regarding integrals around this contour.
We then argue in Section\ \ref{sec:practical} that answering this question can be done with high fidelity by examining integrals of appropriately constructed random functions.
In Section\ \ref{sec:binary} we combine this test with binary search to produce an algorithm capable of bounding the radius of convergence both above and below.
Finally in Section\ \ref{sec:results} we show the results of several numerical experiments performed with these methods and demonstrate that the answers can be made increasingly accurate with increasing computational resources, and that useful results may be extracted with even modest resources.
Notably these methods do not require access to the source code of a function, and so can be applied to black-box analytic functions.
Given the ubiquity of analytic functions and contour integration in scientific computing, we expect these methods to be broadly useful in this context.

\section{Convergence and Integration}
\label{sec:integrate}

Equation\ \eqref{eq:series} combined with Cauchy's Residue Theorem indicates that if there are no poles in a contour $\gamma$ then
\begin{equation}
	f(z_0) = \frac{1}{2\pi i}\oint_{\gamma} \frac{f(z)}{z - z_0} dz,
	\label{eq:mvt}
\end{equation}
where $z_0$ is enclosed by $\gamma$.
This is just the mean value theorem for complex functions.

Equation\ \eqref{eq:mvt} is a necessary condition for $\gamma$ to contain no poles, but it is not sufficient.
To see this consider the function $g(z) = (z^2+1)^{-1}$.
Picking $\gamma$ to be a circular contour of radius $|z-z_0|=2$ around $z_0=0$ gives
\begin{equation}
	\frac{1}{2\pi i}\oint_{\gamma} \frac{g(z)}{z} dz = 1 = g(0),
\end{equation}
so this function satisfies equation\ \eqref{eq:mvt} despite clearly having singularities $z = \pm i$ inside $\gamma$.
This is because the residues of the singularities precisely cancel, so equation\ \eqref{eq:mvt} is insensitive to their presence.

To avoid such cancellation we can multiply $f(z)$ by an analytic function which breaks the cancellation.
In the previous example we could have picked $h(z) = z^2$ and tested $g(z)h(z)$ in addition to $g(z)$.
More generally let $H(\gamma)$ be the set of functions analytic on $\gamma$.
A sufficient condition for $f(z)$ to be analytic on $\gamma$ is that $f(z_0)$ be finite and
\begin{equation}
	f(z_0)h(z_0) = \frac{1}{2\pi i}\oint_{\gamma} \frac{f(z)h(z)}{z - z_0} dz \forall h \in H(\gamma).
	\label{eq:mvt2}
\end{equation}
To see this suppose that $f(z)$ has $N$ poles $z_i$ of order $n_i$ enclosed by $\gamma$.
We take $n_i = 1$ for essential singularities.
Then the functions
\begin{equation}
	h_i(z) = (z - z_0)(z-z_i)^{n_i - 1}\prod_{j\neq i}(z - z_j)^{n_j}
\end{equation}
have the property that
\begin{equation}
	\frac{1}{2\pi i}\oint_{\gamma} \frac{f(z)h_i(z)}{z - z_0} dz = R_i \prod_{j\neq i}(z_i - z_j)^{n_j},
\end{equation}
where $R_i$ is the negative coefficient in the Laurent series for $f(z)$ around this pole.
Noting that $h(z_0)$ is zero and $f(z_0)$ is finite by assumption we see that this is clearly not equal to $f(z_0) h(z_0) = 0$, as each $z_i$ appears precisely once and $R_i \neq 0$.
So long as the number of poles is finite it is clear that $h(z)$ is as well and by construction therefore $h \in H(\gamma)$.
It follows that if there are any poles, and so long as the number of poles is finite, there exists an $h \in H(\gamma)$ such that equation\ \eqref{eq:mvt2} fails to hold.
Thus necessary and sufficient conditions for $f(z)$ to be analytic everywhere in $\gamma$ are:
\begin{enumerate}
\item $f(z_0)$ is finite.
\item $f(z)$ has at most a finite number of poles in $\gamma$.
\item Equation\ \eqref{eq:mvt2} holds for all $h \in H(\gamma)$.
\end{enumerate}

\section{A Practical Test}
\label{sec:practical}

This is an interesting result, but it is far from practical, as it requires examining an uncountable number of contour integrals.
To turn this into a practical test we propose sampling a finite number of random functions from $H(\gamma)$.
In particular, we propose sampling with the functions
\begin{equation}
	h_k (z) = e^{i k (z - z_0)}.
\end{equation}
These are convenient because they are entire (i.e. analytic everywhere), and so they may be used with any contour.
Furthermore they produce non-zero residues even for high-order poles, and they themselves have no roots.

With these functions the integral in equation\ \eqref{eq:mvt2} becomes
\begin{equation}
	\frac{1}{2\pi i}\oint_{\gamma} \frac{f(z)e^{ik(z-z_0)}}{z - z_0} dz = f(z_0) + \sum_j e^{ik (z_j-z_0)}R_j.
\end{equation}
The deviation of this from the expected result is
\begin{equation}
	\frac{1}{2\pi i}\oint_{\gamma} \frac{f(z)e^{ik(z-z_0)}}{z - z_0} dz - f(z_0) h(z_0) = \sum_j e^{ik (z_j-z_0)}R_j.
	\label{eq:mvt3}
\end{equation}
Writing
\begin{equation}
	\boldsymbol{A}(k) \equiv \left\{e^{ik (z_1 - z_0)}, e^{ik (z_2-z_0)}, ..., \right\}
\end{equation}
and
\begin{equation}
	\boldsymbol{R} \equiv \left\{R_1, R_2, ..., \right\}
\end{equation}
we see that the extent to which poles cause a violation of equation\ \eqref{eq:mvt2} is just $\boldsymbol{A}\cdot\boldsymbol{R}$.

Without a prior space on the functions which might be input it is impossible to produce rigorous statistical results regarding the number of $k$'s which will be tested before a violation is found.
Despite this, the process is conceptually intuitive and a rough analysis based on this intuition is possible.
Picking $k$ at random selects a random vector $\boldsymbol{A}(k)$.
While the distribution of this vector is generally not uniform, as poles with $\Im(k z_j) < 0$ are given more weight than those with $\Im(k z_j) > 0$, this effect may be minimised by requiring that
\begin{equation}
	|k| \leq \max_\gamma |z-z_0|.
\end{equation}
Selecting $k$ uniformly subject to this constraint selects $\boldsymbol{A}(k)$ approximately uniformly subject to the constraint that $e^{-1} \leq |A_i| \leq e$.
If the residues $R_j$ are uncorrelated with the pole locations $z_j$ then with $N$ poles the chance that $|\boldsymbol{A} \cdot \boldsymbol{R}| > \epsilon |R| $ for some $\epsilon > 0$ is approximately the chance that two unit vectors on an $N-1$-sphere are not perpendicular to within tolerance $\epsilon$.
This is just
\begin{equation}
	P_\epsilon \approx 1 - \frac{S_{N-2}}{S_{N-1}} \epsilon \approx 1 - \frac{2\sqrt{\pi}}{2N-3}\epsilon,
	\label{eq:prob}
\end{equation}
where $S_{N}$ is the surface area of a unit $N$-sphere.

For large $N$ this method is very good even with modest $\epsilon$.
In fact this allows us to do away with the requirement that our function have a finite number of poles, as equation\ \eqref{eq:mvt3} will almost surely show deviations for one of the sampled functions in this case.
In the case where $N=1$ equation\ \eqref{eq:prob} breaks down, but all such poles are detected with even a single integration because there is no possibility of cancellation.
Somewhat unintuitively then the case of $N=2$ is the most difficult, and there the probability is given by
\begin{equation}
	P_\epsilon \approx 1 - \frac{1}{2\pi} \epsilon.
\end{equation}
for small $\epsilon$.

This is a threshold test, in the sense that we pick a threshold
\begin{equation}
	T \equiv \epsilon |R|
\end{equation}
and compare the integration results against it, reporting pole detection if the integral exceeds $T$ in magnitude and no detection otherwise.
The threshold controls the tradeoff between false-positive and false-negative rates.
False positives, in which we claim a pole detection when there are actually none, may be made less likely by increasing $T$.
False negatives, in which there are poles and we fail to detect them, may be made less likely by decreasing $T$.
In addition, false negatives may be made less likely by performing the test with multiple different $k$'s.
Roughly speaking,
\begin{equation}
	P_{\mathrm{F}-} \approx (1-P_\epsilon)^M \approx \left(\frac{\epsilon}{2\pi}\right)^M
\end{equation}
is the false negative rate, where $M$ is the number of tests performed, and
\begin{equation}
	P_{\mathrm{F}+} \approx M P(T|\delta)
\end{equation}
is the false positive rate, where $P(T|0)$ is the chance that the numerical integration procedure with precision $\delta$ returns a value with magnitude greater than $T$ given that it ought to return zero.

There are then two parameters which can be tuned to minimize these error rates; the number of samples $M$ and the integration accuracy $\delta$.
The tradeoff between the two kinds of errors is controlled by the threshold $T$.
In principle then we can achieve arbitrarily low rates of both errors simply by using more samples and more accurate integration.
The main challenge is that we do not know $|R|$ ahead of time, and so when we select a threshold $T$ we do not know the corresponding $\epsilon$, yet both are relevant parameters for determining the error rate.
To remedy this we approximate
\begin{equation}
	|R| \approx \frac{1}{L(\gamma)}\left|\oint_\gamma |f(z)| dz\right|,
	\label{eq:heurist}
\end{equation}
where $L(\gamma)$ is the length of the contour.
This is inherently heuristic, and indeed can be made to give poor results by such transformations as $f(z) \rightarrow f(z) + 10^{10}$, but most heuristics of this sort suffer from similar issues.
For instance, comparing against the variance of the function along the integration contour is robust against this transformation but fails for the related transformation $f(z) \rightarrow f(z) + 10^{10} z$.
We therefore proceed with this simple heuristic and keep in mind that it is not without shortcomings.

To summarize, this test indicates either that a function is analytic in $\gamma$ or that it is not.
False positives arise due to coincidence in the pole locations and residues with the chosen sample functions, and can be made arbitrarily uncommon by increasing the number of test functions.
False negatives arise only due to insufficiently accurate integration, and so may be made rare by taking $\delta \ll T$.
The only constraints we must impose on $f(z)$ to apply this method are that it can be evaluated everywhere in and on $\gamma$, that it be finite at some point $z_0$ enclosed by $\gamma$, and that it be numerically integrable on $\gamma$.

\section{Search}
\label{sec:search}

Using the test described in the previous section we can determine whether or not a function is analytic in a given contour.
The radius of convergence is the maximum radius passing this test.
Evaluating this radius, however, is generally not practical numerically because a contour with radius equal to the radius of convergence passes through a pole and therefore is quite difficult to evaluate numerically.
This is a fairly generic difficulty in locating poles: if the kind of pole is not known in advance then it is necessary to evaluate the function and its integrals quite near to the pole to properly locate it, and this rapidly becomes an ill-conditioned problem as the desired precision increases.

The next best thing is then to place bounds on the radius of convergence.
In many contexts this is sufficient so long as the bounds are sufficiently tight.
For instance in the case of differentiation by contour integration it typically suffices to place a lower bound within a factor of a few of the radius of convergence\ \citep{Fornberg:1981:NDA:355972.355979}.
Similarly, if the purpose is to identify the nearest pole to a given point then placing bounds on the distance to the pole at least restricts the search space.

A natural algorithm to use in this case is binary search, adapted to the fact that the initial search space is unbounded.
The only complication is that our test has a non-zero error rate in both directions.
One possibility is to use a search algorithm which is tolerant of erroneous answers.
Versions of binary search are known with this property\ \citep{RIVEST1980396, CICALESE2002877}.
These are typically considerably more complicated, however, and come with significant performance overhead.

Given this we choose to simply require that $\delta$ and $\epsilon$ be chosen such that the rates of errors of both kinds are acceptably small.
This incurs a similar performance cost to using more sophisticated search schemes, but has the advantage of being much simpler.
This algorithm is shown in Algorithm\ \ref{algo:search}.
We begin with a trial interval $[0,x]$, where $x$ is a uniform random number drawn from $[0.5,1.5]$ and hypothesize that it contains the radius of convergence.
A random initial upper bound is chosen to make coincidence where the initial contour intersects a pole less likely.
We then expand the upper bound until a contour with that radius contains at least one pole or until some preset limit is reached.
In each subsequent iteration we test the contour at the midpoint of the interval.
If it contains poles we shrink the upper bound of the interval to the contour radius.
If it does not contain poles we set the lower bound of the interval to this radius.
If the integration cannot satisfy the constraint that $\delta \ll T$ with a specified number of samples then we halt the procedure and return the current result.
This result usually means that the test contour is close to the pole.

\begin{algorithm}[H]
\begin{algorithmic}
\State $i \gets 0$
\State $j \gets 0$
\State $L_0 \gets 0$
\State $U_0 \gets \mathrm{rand}()$
\While{$U_0 < \mathrm{Limit}$}
	\If{No poles in circle of radius $U_0$}
		\State $U_{0} \gets 2\times U_{0}$
		\If{$U_{0} > \mathrm{Limit}$}
			\Return 0, $\infty$
		\EndIf
	\Else
		\State $\mathrm{Limit} \gets 0$	
	\EndIf
\EndWhile
\While{$i = 0$ or $U_i - L_j < \mathrm{tolerance}$}
	\State $r \gets (U_i + L_j)/2$
	\If{No poles in circle of radius $r$}
		\State $L_{j+1} \gets r$
		\State $j \gets j+1$
	\ElsIf{Poles in circle of radius $r$}
		\State $U_{i+1} \gets r$
		\State $i \gets i+1$
	\Else{}
		\Return $L_j$, $U_i$
	\EndIf
\EndWhile
\Return $L_j$, $U_i$
\end{algorithmic}
\caption{A binary search procedure.\\Takes as input an upper bound on the radius of convergence ($\mathrm{Limit}$).}
\label{algo:search}
\end{algorithm}

\section{Numerical Results}
\label{sec:results}

For these numerical experiments we perform integration using the trapezoid rule with $N$ points and $M$ samples.
The error in each integral is estimated as the difference between the value obtained with the even-numbered points and that obtained with the odd-numbered points.
This is taken as an estimate for $\delta$, and the procedure halts if $\delta > 0.1 T$, where $T = \epsilon |R|$ and $|R|$ is estimated with equation\ \eqref{eq:heurist}.
The limit on the upper bound is set to $1024$.

Table\ \ref{tab:uniform} shows the results of this method applied to various functions with $\epsilon = 10^{-2}$, $N=10^3$, and $M=3$.
Our method succeeds even for functions such as $(1+z)^{-1} + (1-z)^{-1}$, which an unmodified contour integral would erroneously conclude to be pole-free for $r > 1$.
For simple poles of low order the algorithm performs quite well, typically bounding the radius of convergence to within a 10\%.

Artificially inserting worst-case terms for the heuristic in equation\ \eqref{eq:heurist} causes this performance to degrade, but overall it is still quite good.
For essential singularities and high-order poles the algorithm still returns correct bounds, but they are typically somewhat wider, with $U/L$ of order $2$.
The fact that no incorrect results were returned indicates that $\epsilon$ and $\delta$ were chosen sufficiently conservatively.

\begin{table}
\begin{tabular}{|c|c|c|c|c|}
\hline 
Function & $z_{0}$ & $r_{\mathrm{convergence}}$ & Lower bound & Upper bound\tabularnewline
\hline
\hline 
$\frac{1}{1+z}$ & $0$ & $1$ & $0.984$ & $1.023$\tabularnewline
\hline 
$\frac{1}{1+z}$ & $2$ & $3$ & $2.913$ & $3.059$\tabularnewline
\hline 
$\frac{1}{1+z}$ & $-1.1$ & $0.1$ & $0.0959$ & $0.102$\tabularnewline
\hline 
$\frac{1}{1+z}$ & $i$ & $\sqrt{2}\approx1.414$ & $1.394$ & $1.416$\tabularnewline
\hline 
$\frac{1}{1+z} + \frac{1}{1-z}$ & $0$ & $1$ & $0.806$ & $1.209$\tabularnewline
\hline 
$10^2 + \frac{1}{1+z}$ & $0$ & $1$ & $0.923$ & $1.077$\tabularnewline
\hline 
$\frac{1}{1+z}+\frac{2}{1-z^{2}}$ & $0$ & $1/\sqrt{2}\approx0.707$ & $0.683$ & $0.716$\tabularnewline
\hline 
$\exp\left[\frac{1}{1+z}\right]$ & $0$ & $1$ & $0.900$ & $1.050$\tabularnewline
\hline 
$z$ & $0$ & $\infty$ & $0$ & $\infty$\tabularnewline
\hline 
$\cos\left(\frac{1}{z}\right)$ & $i$ & $1$ & $0.659$ & $1.318$\tabularnewline
\hline 
$\frac{z^{8}}{z^{12}+1}$ & $4$ & $5$ & $2.877$ & $5.755$\tabularnewline
\hline
\end{tabular}
\label{tab:uniform}
\caption{Results of the algorithm applied to various test functions with $M=3$ samples, $N=10^3$ points per integration, and $\epsilon=10^{-2}$.}
\end{table}

\begin{figure}
	\centering
	\includegraphics[width=0.9\textwidth]{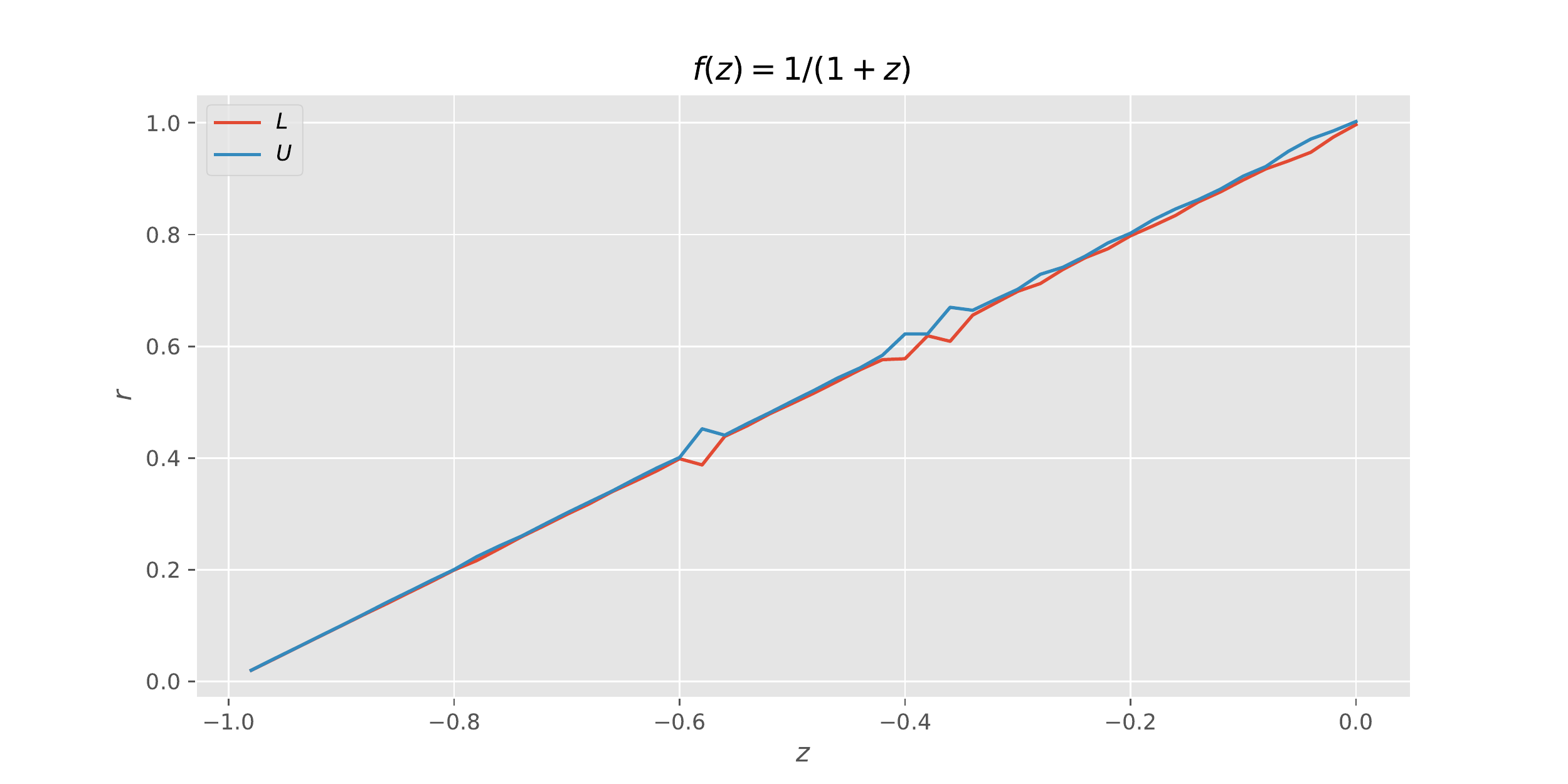}
	\caption{The upper and lower bounds on the radius of convergence is shown as a function of the reference position $z_0$. For this test $N=10^4$ and $M=3$.}
	\label{fig:z}
\end{figure}

Figure\ \ref{fig:z} shows the upper and lower bounds on the radius of convergence as a function of the reference position $z_0$ for the function $(1+z)^{-1}$.
The true radius of convergence traces a straight line in between the two curves.
There are small but sharp variations in the tightness of the bounds.
These occur when the search happens to terminate a step earlier or later than the searches on either side of it, leading to a difference of a factor of $2$ between neighboring points.

\begin{figure}
	\centering
	\includegraphics[width=0.9\textwidth]{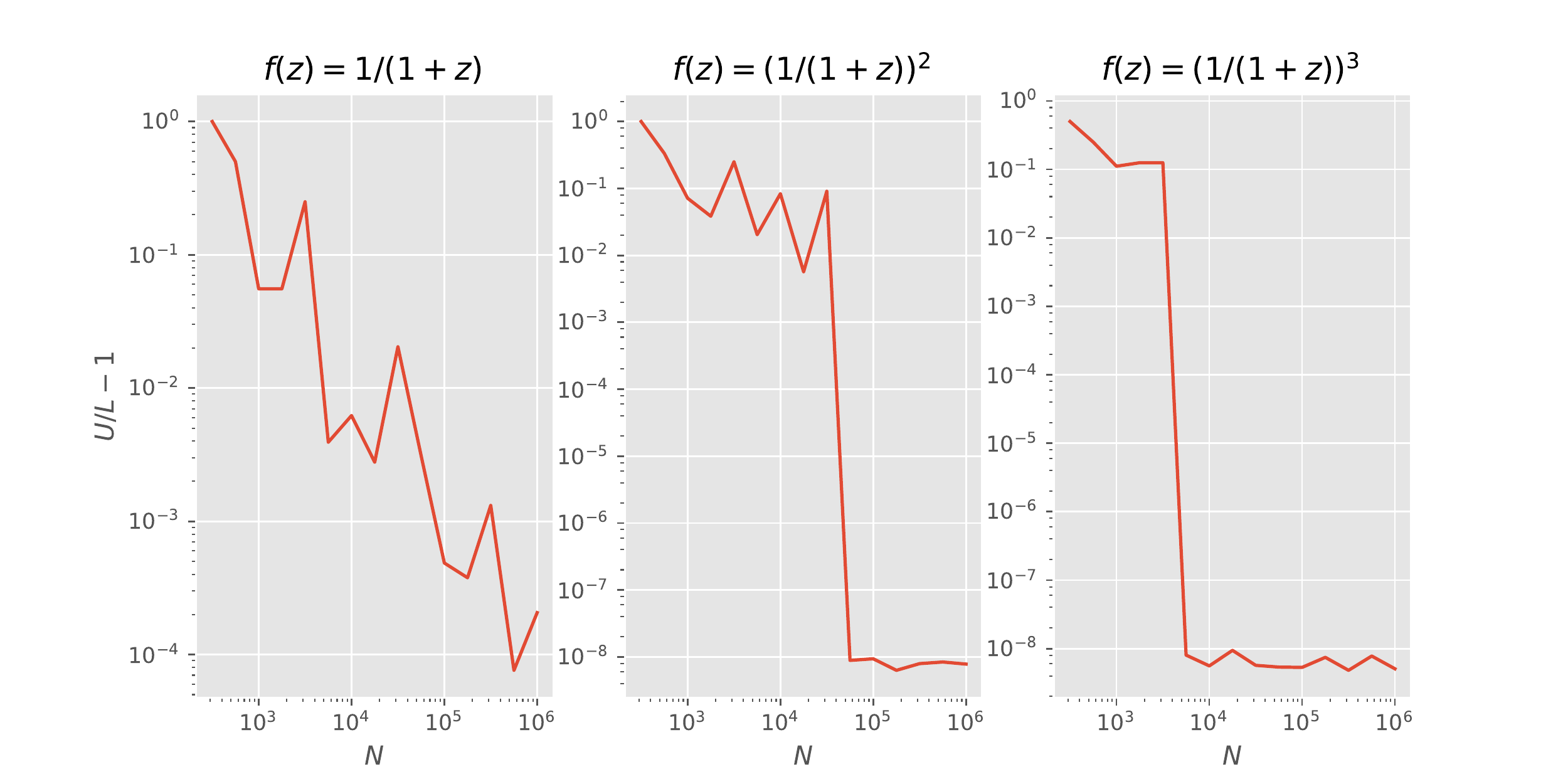}
	\caption{The tightness of the bounds is shown as a function of $N$ on a log-log scale.}
	\label{fig:N}
\end{figure}

Figure\ \ref{fig:N} shows $U/L - 1$ versus $N$ for several functions with $M=5$.
There is a strong downward trend, indicating that the bounds become quite tight as $N$ increases.
The precise shape of this trend depends on the order of the pole.
For the first-order pole the bound scales as $N^{-1}$, which is consistent with the fact that the integral is dominated by a region near the pole of width proportional to the distance to the pole.
For the essential singularities the precise scaling is less clear, but the bounds still become tight at large $N$.

\section{Conclusions}

We have described a novel method for robustly identifying whether or not a given function has any poles inside a given contour.
We used this to develop a method for placing both upper and lower bounds on the radius of convergence of complex-analytic functions.
This method is well-behaved in numerical experiments, producing tight bounds for simple poles and accurate, if less tight, bounds for tighter poles and essential singularities.
This is the key utility of our method: it can be applied to black-box functions without access to their symbolic representations, and it produces accurate and useful results even when given functions with mild pathologies.
This algorithm may therefore be used broadly and with little tuning in problems which require interrogating the analytic structure of functions.

\section*{Acknowledgments}
The author gratefully acknowledges helpful comments from Ravishankar Sundararaman as well as financial support from the United Kingdom Marshall Commission. 

\bibliography{refs}

\end{document}